\documentclass{amsart}
\usepackage{amsmath,amssymb}
\usepackage[all,dvips]{xy}
\newcommand{\PP}{\mathbb{P}}
\newcommand{\ZZ}{\mathbb{Z}}
\newcommand{\QQ}{\mathbb{Q}}

\newcommand{\Aut}{\mathrm{Aut}}
\newcommand{\Gal}{\mathrm{Gal}}

\newcommand{\Jac}{\mathrm{Jac}}

\newcommand{\genus}{\mathrm{genus}}

\newcommand{\res}{\mathrm{res}}
\newcommand{\tors}{\mathrm{tors}}
\newcommand{\rk}{\mathrm{rk}}
\newcommand{\Kbar}{{\overline{K}}}

\newcommand{\act}[2]{{}^{#1}\!#2}

\newcommand{\sM}{\mathcal{M}}

\newlength{\underscorelength}
\newcommand{\deadus}{\settowidth{\underscorelength}{{\rm \_}}\hspace{-0.5\underscorelength}\_\hspace{-0.5\underscorelength}}
\newcommand{\Sha}{{\rm I{\deadus}I{\deadus}I}}

\renewcommand{\phi}{\varphi}
\newtheorem{theorem}{Theorem}[section]
\newtheorem{lemma}[theorem]{Lemma}
\newtheorem{corollary}[theorem]{Corollary}
\newtheorem{prop}[theorem]{Proposition}
\theoremstyle{definition}

\begin{document}
\title{Visualising Sha[2] in Abelian Surfaces}
\parskip=0.1pt
\author{Nils Bruin}
\address{School of Mathematics, University of Sydney NSW 2006, Australia}
\email{bruin@member.ams.org}
\thanks{The research in this paper was funded by the Pacific Institute for the
Mathematical Sciences, Simon Fraser University, the University of British
Columbia and the School of Mathematics of the University of Sydney.
}


\date{September 10, 2002}

\begin{abstract}
Given an elliptic curve $E_1$ over a number field and an element $s$ in its $2$-Selmer
group, we give two different ways to construct infinitely many Abelian surfaces
$A$ such that the
homogeneous space representing $s$ occurs as a fibre of $A$ over another
elliptic curve $E_2$.
We show that by comparing the $2$-Selmer groups of $E_1$, $E_2$ and $A$, we can
obtain information about $\Sha(E_1/K)[2]$ and we give examples where we use this
to obtain a sharp bound on the Mordell-Weil rank of an elliptic curve.

As a tool, we give a precise description of the $m$-Selmer group of an Abelian
surface $A$ that is $m$-isogenous to a product of elliptic curves $E_1\times
E_2$.

One of the constructions can be applied iteratively to obtain information about
$\Sha(E_1/K)[2^n]$. We give an example where we use this iterated application
to exhibit an element
of order $4$ in $\Sha(E_1/\QQ)$.
\end{abstract}
\maketitle
\section{Introduction}

The Mordell-Weil theorem states that the rational points on an elliptic curve
$E$ over a number field $K$ form a finitely generated commutative group. Given
$m\in\ZZ$ with $m\geq 2$, there is an in principle effectively computable
object, the $m$-Selmer group $S^{(m)}(E/K)$, that provides an upper bound on
the free rank of $E(K)$, the \emph{$m$-Selmer-rank} of $E/K$. This bound needs
not be sharp. The $m$-torsion of the \emph{Shafarevich-Tate} group,
$\Sha(E/K)[m]$, measures the failure of the $m$-Selmer rank to provide a sharp
bound on the rank of $E(K)$. If $\Sha(E/K)[m]$ has no elements of order $m$,
then the Selmer-rank equals the rank of $E(K)$.

Suppose $E_1$ and $E_2$ are elliptic curves over a number field $K$ with
$E_1[m]\simeq E_2[m]$. We write $\Delta:=E_1[m]\simeq E_2[m]$.
One can construct an Abelian surface $A=E_1\times
E_2/\Delta_E$, where $\Delta_E\subset\Delta\times\Delta\subset E_1\times E_2$ is
the anti-diagonal in $\Delta\times\Delta\simeq E_1[m]\times E_2[m]$.
We investigate how
$\Sha(E_1/K)[m]$ and $\Sha(E_2/K)[m]$ are related to $\Sha(A/K)[m]$.
In particular, we prove

\begin{theorem}
\label{thrm:biell}
Let $E_1$ be an elliptic curve over a number field $K$. Let
$\xi\in
\Sha(E/K)[2]$. There are infinitely many explicitly constructible elliptic
curves $E_2$, not isomorphic over $\Kbar$, such that $\xi$ is in the kernel of the
natural map $\Sha(E_1/K)[2]\to \Sha(A/K)[2]$.
\end{theorem}

A non-constructive proof of the existence of $A$ can be found in
\cite{klenke:visWC2}. In fact, the proof there applies to any $\xi\in
H^1(K,E[2])$.

One obvious choice is to take $E_2$ to be a quadratic twist of $E_1$. Suppose
that $E_2$ is the twist of $E_1$ by $d\in K^*$. Then $E_2$ is isomorphic to
$E_1$ over $K(\sqrt{d})$ the Abelian surface $A$ is the Weil-restriction of
$E_1$ with respect to $K(\sqrt{d})/K$.

\begin{theorem}
\label{thrm:quadext}
Let $E$ be an elliptic curve over a number field $K$. Let $\xi\in
\Sha(E/K)[2]$. There are infinitely many explicitly constructible quadratic
extensions $K(\sqrt{d})/K$ such that $\xi$ is in the kernel of
$\Sha(E/K)[2]\to\Sha(E/K(\sqrt{d}))[2]$.
\end{theorem}

A non-constructive proof that any $\xi\in\Sha(E/\QQ)[m]$, with $m\geq 2$, is
visualised in an Abelian variety of dimension at most $m$ can be found in
\cite{steinaga:visabvar}. We can use Theorem~\ref{thrm:quadext} to get an
explicit version of this fact for $m=2^n$ and arbitrary base field.

\begin{corollary}
\label{cor:twopowvis}
Let $E$ be an elliptic curve over a number field $K$. Let $\xi\in
\Sha(E/K)[2^n]$. Then there is an explicitly constructible $n$-dimensional
Abelian variety $A$ over $K$
with a non-constant map $E\to A$ such that $\xi$ vanishes under the natural
map $H^1(K,E)\to H^1(K,A)$.
\end{corollary}

\begin{proof}

Let $K_0=K$ and $\xi_0=\xi$. By Theorem~\ref{thrm:quadext}, there is a
quadratic extension $K_1$ of $K$ such that $2^{(n-1)}\xi_0\in\Sha(E/K_0)[2]$
vanishes in $\Sha(E/K_1)$. Consequently, the image $\xi_1$ of $\xi_0$ in
$\Sha(E/K_1)$ satisfies $\xi_1\in\Sha(E/K_1)[2^n-1]$. We repeat this
construction for $i=0,1,\ldots,$ to obtain a quadratic extension $K_{i}$ of
$K_{i-1}$ such that the image $\xi_i$ of $\xi_{i-1}$ in $\Sha(E/K_{i})$
satisfies $2^{n-i}\xi_{i}=0$. It follows that $\xi_n=0$ in $\Sha(E/K_n)$.

Let $A=\Re_{K_n/K} E$ be the Weil-restriction of scalars in the sense of
\cite[\S 7.6]{BLR:neron}. Then $A$ is a $2^n$-dimensional Abelian variety and
there is a non-constant morphism $E\to A$. By Shapiro's Lemma
(\cite[Proposition~2]{atwall:groupcohom}) we have, canonically,
$H^1(K_n,E)\simeq H^1(K,A)$, and the image of $\xi$ under $H^1(K,E)\to H^1(K,A)$
is indeed trivial.
\end{proof}

For the proofs of Theorems~\ref{thrm:biell} and \ref{thrm:quadext}, we analyse
the cohomology of non-simple Abelian surfaces.
In particular, let $K$ be a number field and let
$A$ be Abelian surface over $K$ with a degree $m$ isogeny
$\phi:A\to E_1\times E_2$ to a product of elliptic curves $E_1$ and $E_2$. We
show that the Selmer groups satisfy the following equalities.
$$\begin{array}{rcl}
S^{(\phi)}(A/K)&=&S^{(m)}(E_1/K)+S^{(m)}(E_2/K)\\
S^{(\hat{\phi})}(E_1\times E_2/K)&=&S^{(m)}(E_1/K)\cap S^{(m)}(E_2/K)
\end{array}$$

As a motivation and an illustration, we consider a simple example that
illustrates both theorems and shows how the constructions can be applied to
exhibit non-trivial elements in $\Sha(E/K)[2]$.

Consider the elliptic curve over $\QQ$ given by
$$E_1:y^2=x^3-22x^2+21x+1.$$
From the $2$-Selmer group of $E_1/\QQ$, we see that
$\rk E_1(\QQ)\leq 4$. A simple computation
shows that $\langle(0,1),(1,1)\rangle\subset E_1(\QQ)$ forms a free subgroup of
rank $2$. An an analytic rank computation and a $3$-descent (see
\cite{schasto:pdescent}) yield convincing evidence that the rank of $E_1$ really
is $2$, which means that $\Sha(E_1/\QQ)[2]\simeq\ZZ/2\ZZ\times\ZZ/2\ZZ$.

\begin{prop}\label{prop:schasto}
$\rk E_1(\QQ)=2$.
\end{prop}

\noindent\emph{Proof 1}: Consider
$$E_2: 2y_2^2=x^3-22x^2+21x+1.$$
The group $\langle (1/2,7/4),(1/8,41/32)\rangle\subset E_2(\QQ)$ is free of rank
$2$, so $\rk E_2(\QQ)\geq 2$. From the $2$-Selmer group of $E_2/\QQ$ we deduce
that $\rk E_2(\QQ)\leq 4$.
 From
$$\rk E_1(\QQ(\sqrt{2}))=\rk E_1(\QQ)+\rk E_2(\QQ),$$
it follows that $\rk E_1(\QQ(\sqrt{2}))\geq 4$. From the $2$-Selmer group of
$E_1/\QQ(\sqrt{2})$, we also find $4$ as an upper bound for $\rk
E_1(\QQ(\sqrt{2}))$. Consequently, we have $\rk E_1(\QQ)=\rk E_2(\QQ)=2$.
\hfill$\square$
\medskip

\noindent\emph{Proof 2}: Consider the smooth complete curve corresponding to
the affine model
$$C: y_1^2=-y_0^6-19y_0^4+20y_0^2+1.$$
Let
$$E_2: v^2=u^3+20u^2-19u-1.$$
The curve $C$ of genus $2$ covers $E_1$ by
$(y_0,y_1)\mapsto(x,y_1)=(-y_0^2+1,y_1)$ and $E_2$ by
$(y_0,y_1)\mapsto(u,v)=(1/y_0^2,y_1/y_0^3)$. It follows that $\Jac_C$ is
isogenous to $E_1\times E_2$ and that
$$\rk \Jac_C(\QQ)=\rk E_1(\QQ)+\rk E_2(\QQ).$$
From $\langle (1,1),(2,7),(5,23) \rangle\subset E_2(\QQ)$ and the $2$-Selmer
group of $E_2/\QQ$, we deduce that $\rk E_2(\QQ)=3$. A $2$-descent on
$\Jac_C/\QQ$ yields $\rk \Jac_C(\QQ)=5$, so it follows that $\rk E_1(\QQ)=2$.
\hfill$\square$
\medskip

Both proofs make use of essentially the same construction. We find an elliptic
curve $E_2$
and an Abelian surface $A$ isogenous to $E_1 \times E_2$. Using a $2$-descent we
show that the rank of $A(\QQ)$ is smaller than the sum of the rank-bounds we get
from a $2$-descent on $E_1$ and $E_2$ separately. In this article we analyse
when this construction may yield non-trivial results and we give an explicit
construction for $E_2$.

\section{Selmer groups}
\subsection{Abstract definition of the Selmer group}
In this section we recall the abstract definition of the Selmer group
associated to an isogeny between Abelian varieties. We review the relation of
the size of Selmer groups to the rank of an Abelian variety over a number field.
Given isogenies $\phi$ and $\hat{\phi}$ such that $\phi\circ\hat{\phi}$ is
multiplication-by-$m$,
we indicate how one can use the full multiplication-by-$m$ Selmer group to
improve the rank-bound obtained from the $\phi$- and $\hat{\phi}$-Selmer groups.

Consider two Abelian varieties $A$ and $B$ of equal dimension over a field $K$
of characteristic
$0$ and a finite morphism (an \emph{isogeny}) $\phi:A\to B$. Let
$\Delta=\ker\phi$. Galois-cohomology yields
$$0\to B(K)/\phi A(K)\to H^1(K,\Delta)\to H^1(K,A).$$

If $K$ is a number field, then we can approximate
the image of $B(K)/\phi A(K)$ in $H^1(K,\Delta)$ by determining
the elements of $H^1(K,\Delta)$ that are in the image everywhere
locally. This constitutes the $\phi$-Selmer group of $A$ over $K$. The following
diagram with exact rows illustrates the definition. The products are
taken over all primes $p$ of $K$.
$$\xymatrix{
0 \ar[r] &
  B(K)/\phi A(K) \ar[r] \ar[d] &
  H^1(K,\Delta) \ar[r] \ar[d] &
  H^1(K,A) \ar[d] \\
0 \ar[r] &
  \prod_p B(K_p)/\phi A(K_p) \ar[r] &
  \prod_p H^1(K_p,\Delta) \ar[r] &
  \prod_p H^1(K_p,A)
}$$
The $\phi$-Selmer group of $A$ over $K$ is defined to be the subgroup of
$H^1(K,\Delta)$ consisting of cocycles that map to elements in
$\prod_p H^1(K_p,\Delta)$ that have a pre-image in $\prod_p B(K_p)/\phi A(K_p)$.
It is defined by the exact sequence
$$0\to S^{(\phi)}(A/K)\to H^1(K,\Delta)\to \prod_p H^1(K_p,A).$$
The Selmer-group contains $B(K)/\phi A(K)$ and therefore provides a bound on its
size. Unfortunately, there may be $1$-cocycles of $A$ that are trivial
everywhere locally (restrict to a coboundary under $H^1(K,A)\to H^1(K_p,A)$ for
all primes $p$ of $K$), but are not coboundaries themselves. They form the
\emph{Shafarevich-Tate group}
$$0\to \Sha(A/K)\to H^1(K,A)\to \prod_p H^1(K_p,A).$$

The subgroup of everywhere locally trivial cocycle classes that are in the
kernel of $H^1(K,A)\to H^1(K,B)$ measures the failure of $S^{(\phi)}(A/K)$ to
bound $B(K)/\phi A(K)$ sharply:
$$0\to B(K)/\phi A(K)\to S^{(\phi)}(A/K)\to\Sha(A/K)[\phi]\to 0.$$
Although from the description given here, it is not clear that either
$S^{(\phi)}(A/K)$ or $\Sha(A/K)[\phi]$ are effectively computable, we will
shortly see that the Selmer group is finite and effectively computable in the
situations we are interested in. First we explain why $\phi$-Selmer groups help
in computing the rank of $A(K)$.

First suppose that $\phi$ is multiplication by $m$, so that $B=A$.
By The Mordell-Weil theorem we know that $A(K)\simeq \ZZ^r\times A(K)_\tors$,
where $A(K)_\tors\subset A(K)$ is the finite subgroup of elements of finite
order. Consequently,
$$A(K)/mA(K)=(\ZZ/m\ZZ)^r\times A(K)_\tors/mA(K)_\tors.$$
Since the group $A(K)_\tors$ is generally relatively easy to compute, one can
deduce the free rank $r$ from the size of $A(K)/mA(K)$.

If $m$ is prime, we have
$\#A(K)_\tors/mA(K)_\tors=\#A[m](K)$. For an isogeny $\phi:A\to B$ such that
$\hat{\phi}\phi=m$ with $m$ prime, we define, analogous to the fact that $m^{\rk
A(K)}\#A[m](K)=\#A(K)/mA(K)$, the number $s=:\rk S^{(\phi)}(A/K)$ by
$m^s\#A[\phi](K)=\#S^{(\phi)}(A/K)$. This definition allows us to state the
relation between the rank of an Abelian variety and its $m$-Selmer group
concisely.

\begin{prop} Let $A$ be an Abelian variety over a number field $K$ and $p$ a
prime number. Then
$\rk A(K)\leq \rk S^{(p)}(A/K)$ and equality holds precisely if
$\#\Sha(A/K)[p]=1$.
\end{prop}

Selmer groups of other isogenies also give information about the rank by
combining them with the Selmer group of the dual isogeny.

\begin{lemma}
\label{lemma:cardiso}
Let $\phi:A\to B$ be an isogeny of Abelian varieties over a field
$K$ such that $\#B(K)/\phi A(K)$ and $\#A(K)/\hat{\phi}B(K)$ are finite groups.
Suppose that $\phi\hat{\phi}=m$. Then
$$
\frac{\# A(K)/m A(K)}{\# A[m](K)}=
\frac{\#B(K)/\phi A(K)}{\#A[\phi](K)}
\frac{\#A(K)/\hat{\phi} B(K)}{\#B[\hat{\phi}](K)}.
$$
\end{lemma}
\begin{proof}
Consider the exact sequence of finite groups
$$\begin{array}{l}
0\to A[\phi](K)\to A[m](K)\to B[\hat{\phi}](K)\to\\
\quad B(K)/\phi A(K)\to A(K)/m A(K)\to A(K)/\hat{\phi} B(K)\to 0.
\end{array}$$
\end{proof}

\begin{corollary} \label{cor:dwlsel}
Let $\phi: A\to B$ be an isogeny of Abelian varieties over a number
field $K$ such that $\hat{\phi}\phi=p$ for some prime $p$. Then
$$\rk A(K)\leq \rk S^{(\phi)}(A/K) + \rk S^{(\hat{\phi})}(B/K)$$
and equality holds precisely if $\#\Sha(A/K)[\phi]=\#\Sha(B/K)[\hat{\phi}]=1$.
\end{corollary}

\subsection{Two-Selmer groups of elliptic curves}
\label{sec:selell}

In this section we give two alternative descriptions of $H^1(K,E[2])$ and the
map $E(K)/2E(K)\to H^1(K,E[2])$,
which we will use in Sections~\ref{sec:quadext} and \ref{sec:biell}.

To find $S^{(2)}(E/K)\subset H^1(K,E[2])$, one
needs some further non-trivial computations, see for instance
\cite{stoll:descent}.

We consider an elliptic curve
$$E:y^2=x^3+a_2 x^2+a_4 x+a_6=F(x)$$
over a field $K$ of characteristic $0$. First we describe
$H^1(K,E[2])$.
We write $M_K=K[x]/(F(x))$ and let $\theta\in M$ denote the residue class of
$x$, i.e., $F(\theta)=0$. We write $M_K^*$ for the multiplicative group of the
algebra $M_K$, $M_K^{*2}\subset M_K^*$ for the squares and $M_K'\subset M_K^*$
for the kernel of $M_K^*\stackrel{N_{M_K/K}}{\longrightarrow} K^*/K^{*2}$.

\begin{theorem}[{\cite{cassels:diophell}, \cite[Chapter~15]{cas:ell}}]
\label{thrm:H1E} Let $K$ be a field of characteristic $0$ and let $E$ and
$M_K$ be as above.
Then
$$H^1(K,E[2])\simeq M_K'/M_K^{*2}.$$
Under this isomorphism, the map
$$\mu:E(K)/2E(K)\to H^1(K,E[2])$$
is induced by
$$\begin{array}{cccl}
E(K)&\to&M_K^*\\
(x,y)& \mapsto & x-\theta &\mbox{if $F(x)\neq 0$,}\\
\end{array}$$
\end{theorem}

To facilitate the evaluation of the map $\mu$, we use the following lemma, which
can be proved by a straightforward computation.

\begin{lemma}
\label{lemma:thetasqr}
Let
$E: y^2=x^3+a_2x^2+a_4x+a_6=F(x)$ be an elliptic curve over a field $K$
and let $\theta$ be a root of $F(x)$ in $K[x]/F(x)$.
Let
$e_1,e_2\in E(K)$. Then
$$(x(e_1)-\theta)(x(e_2)-\theta)(x(e_3)-\theta)=C(e_1,e_2)^2$$
where $C(e_1,e_2)$ is a symmetric algebraic function of $e_1,e_2$.
\end{lemma}

\begin{corollary}
\label{cor:munorm}
Let $E: y^2=x^3+a_2x^2+a_4x+a_6=F(x)$ be an elliptic curve over a field $K$ and
let $L$ be a separable quadratic extension of $K$, with $\sigma:L\to L$
conjugation over $K$. Then the following diagram commutes.
$$\xymatrix{
E(L) \ar[r]^\mu \ar[d]^{e\mapsto e+\act{\sigma}{e}} &
  M_L'/M_L^{*2} \ar[d]^{N_{M_L/M_K}}\\
E(K) \ar[r]^\mu & M_K'/M_K^{*2}}$$
\end{corollary}
\begin{proof}
The only possible obstruction to the commutativity of this diagram is that
$M_L^{*2}\cap M_K^*$ may be larger than $M_K^{*2}$. However,
using Lemma~\ref{lemma:thetasqr}, we see that, for $e$ such that
$C(e,\act{\sigma}{e})$ is well defined and non-zero,
$\mu(e+\act{\sigma}{e})\mu(e)\mu(\act{\sigma}{e})=C(e,\act{\sigma}{e})^2$,
where $C(e,\act{\sigma}{e})$ is $\sigma$-invariant and therefore in $M_K^*$
For the special $e$ such that $C(e,\act{\sigma}{e})$ has a pole or a zero, one
can verify the statement separately, which we will leave to the reader.
\end{proof}

An alternative interpretation of $H^1(K,E[2])$ can be obtained by considering
unramified covers of $E$ with a certain Galois-group. First we review some
terminology.
Let $\pi: D\to C$ be a non-constant morphism of smooth complete absolutely 
irreducible
curves over $K$ such that
$\#\Aut_\Kbar(D/C)=\deg\phi$, i.e., a \emph{Galois cover}. We say that $\pi': D'\to C$ is a \emph{twist} of
$D\stackrel{\pi}{\to}C$ if there is an isomorphism $\psi: D\to D'$ over $\Kbar$ such
that $\pi=\pi'\circ\psi$. If $\psi$ is already defined over $K$ then
$D\stackrel{\pi}{\to} C$ and $D'\stackrel{\pi'}{\to} C$ are considered the same.
\begin{theorem}[\cite{sil:AEC1}] Let $\pi:D\to C$ be a Galois cover over a
number field $K$. Then
$$H^1(K,\Aut_\Kbar(D/C))\simeq\{\mbox{Twists of }D\stackrel{\pi}{\to}C\}$$
as pointed sets with $Gal(\Kbar/K)$-action.
\end{theorem}

In particular, $H^1(K,E[2])\simeq \{\mbox{Twists of }E\stackrel{2}{\to}E\}$.
Given an element $\delta\in M_K'$, we can find the corresponding cover of $E$
explicitly. We refer to it as $T_\delta$ to distinguish it from the
representation as an element of $M_K'/M_K^{*2}$. First, we fix a representation of
$M_K$ over $K$. Let $\{1,\theta,\alpha\}$ be a basis of $M_K$ as a $K$-vector space. If
$F(x)$ is irreducible over $K$, then one can take $\alpha=\theta^2$.
We find a model of $T_\delta$ by composing $T_\delta\to E$ with
$E\stackrel{x}{\to}\PP^1$. A heuristic motivation for the construction below is
that if $P\in E(K)$ with $\mu(P)=\delta$, then $T_\delta$ should have a rational
point above $P$, so in the composed cover $T_\delta\to\PP^1$, we have a rational
point above $x(P)$.

If $\mu(P)=\delta$, then there are $u_0,u_1,u_2\in K$ such that
$$x(P)-\theta=\delta(u_0+\theta u_1+\alpha u_2)^2.$$
We expand the right hand side with respect to the $K$-basis of $M_K$.
Let $Q_{\delta,i}(u_0,u_1,u_2)\in K[u_0,u_1,u_2]$
be the quadratic forms so that
$$\delta(u_0+\theta u_1+\alpha u_2)^2=Q_{\delta,0}(u_0,u_1,u_2)+
     \theta Q_{\delta,1}(u_0,u_1,u_2)+
     \alpha Q_{\delta,2}(u_0,u_1,u_2).$$
We find that $Q_{\delta,0}(u_0,u_1,u_2)=x(P)$, $Q_{\delta,1}(u_0,u_1,u_2)=-1$
and $Q_{\delta,2}=0$.     

Independent of whether there exists a point $P\in E(K)$ with $\mu(P)=\delta$,
we take
$u_0,u_1,u_2$ to be variables and define the $Q_{\delta,i}$ as above.
Using these forms, we define two projective varieties over $K$. Let
$$\begin{array}{rl}
L_\delta:&Q_{\delta,2}(u_0,u_1,u_2)=0,\\
T_\delta:&Q_{\delta,2}(u_0,u_1,u_2)=0\mbox{ and }
                               Q_{\delta,1}(u_0,u_1,u_2)=-u_3^2.
\end{array}$$
Note that the map $(u_0:u_1:u_2:u_3)\mapsto(u_0:u_1:u_2)$ induces a degree $2$
cover $T_\delta\to L_\delta$.
This cover is ramified 
ramified at the points where $u_3=0$. Consequently, the ramification locus lies
above those $(u_0:u_1:u_2)$ satisfying
$Q_{\delta,1}(u_0,u_1,u_2)=Q_{\delta,2}(u_0,u_1,u_2)=0.$
On both $L_\delta$ and $T_\delta$ we define the function
$$x(u_0,u_1,u_2)=-\frac{Q_{\delta,0}(u_0,u_1,u_2)}{Q_{\delta,1}(u_0,u_1,u_2)}.$$
Since $\delta\in M_K'$, we can choose $d\in K^*$ such that $d^2=N_{M_K/K}(\delta)$.
On $T_\delta$, we define
$$y_d(u_0,u_1,u_2)=\frac{d}{u_3^3}N_{M_K[u_0,u_1,u_2]/K[u_0,u_1,u_2]}(u_0+\theta u_1+\alpha
u_2)$$
We find that $T_\delta$ covers $E$ by
$$\begin{array}{rccc}
\phi_\delta:&T_\delta&\to&E\\
&(u_0:u_1:u_2:u_3)&\mapsto&(x(u_0,u_1,u_2),y_d(u_0,u_1,u_2,u_3))
\end{array}$$
The defined varieties are in fact connected curves and they fit in the
following diagram.
$$\xymatrix@R-2em{
&T_\delta \ar[dr] \ar[ddl]^\phi\\
&&L_\delta \ar[ddl]_x \\
E \ar[dr]_x \\
& \PP^1
}$$
We see that $T_\delta=E\times_{\PP^1} L_\delta$. One can check that
$T_\delta\to E$ is unramified. Consequently, $\genus(T_\delta)=1$. Furthermore, 
$T_\delta\to L_\delta$ is ramified at $4$ geometric points, so
$\genus(L_\delta)=0$.
We have $\delta\in S^{(2)}(E/K)$ precisely if $T_\delta(K_p)$ is non-empty for
all primes $p$ of $K$. Consequently, $L_\delta(K_p)$ is also non-empty for all
primes $p$ of $K$, and by the Hasse-Minkowsky-theorem, we have $L_\delta\simeq
\PP^1$. Choose a parametrisation  $t\mapsto(u_0(t),u_1(t),u_2(t))$ of
$L_\delta$. This yields a model
$$T_\delta:u_3^2=R(t)=-Q_{\delta,1}(u_0(t),u_1(t),u_2(t)).$$
We recover the description of $S^{(2)}(E/K)$ as a set of classes of quartics
that have a point everywhere locally, as used in some formulations of the
$2$-descent method for elliptic curves over $\QQ$ (\cite{simon:ranks},
\cite{cremona:clasinv}, \cite{cremona}).

\begin{prop}\label{thrm:covrep}
Let $E$, $K$, $M_K'$ and $T_\delta$ be defined as above. The isomorphism
$$H^1(K,E[2])\to\mathrm{Twists}(E\stackrel{2}{\to}E)$$
is induced by
$$\begin{array}{ccc}
M_K'&\to&\{\mbox{Covers of } E\}\\
\delta&\mapsto&T_\delta.
\end{array}$$
\end{prop}

\section{Non-simple Abelian surfaces}
\label{sec:cohom}

In this section we consider an Abelian surface $A$, isogenous to the product of
two elliptic curves $E_1$ and $E_2$. Although in Sections~\ref{sec:quadext} and
\ref{sec:biell} we will only be considering $2$-isogenies, the results in this
section are valid for any $m$-isogeny $m\geq 2$.

Let $m>1$ be an integer and let $E_1$ and $E_2$ be elliptic curves over a
field $K$ of characteristic $0$ with $E_1[m]\simeq E_2[m]$ as group schemes
over $K$. We write
$\Delta:=E_1[m]\simeq E_2[m]$ and $\Delta_E\subset E_1[m]\times E_2[m]$
for the anti-diagonal embedding. Let $A:=(E_1\times E_2)/\Delta_E$.

We have $p^*:E_1\to A$ induced by $P\mapsto (P,0)\in E_1\times E_2$ and
$q^*:E_2\to A$ induced by $Q\mapsto (0,Q)\in E_1\times E_2$. 
Note that $\Delta\subset (E_1\times E_2)[m]$, so the multiplication-by-$m$
map on $E_1\times E_2$ factors through $A$. This factorisation induces the maps
$p_*:A\to E_1\times E_2\to E_1$ and $q_*: A\to E_1\times E_2\to E_2$. It is
straightforward to check that $p_*\circ p^*=m|_{E_1}$ and $q_*\circ
q^*=m|_{E_2}$.
$$\xymatrix{
&&&0 \ar[d] \\
&&&\Delta_A \ar[d] \\
0 \ar[r] &
  \Delta_E \ar[r] &
  E_1 \times E_2 \ar[r]^{p^*+q^*} \ar[dr]_m &
  A \ar[r] \ar[d]^{p_*\times q_*} &
  0\\
&&& E_1\times E_2 \ar[d] \\
&&& 0
}$$
It follows that the isogeny $p_*\times
q_*:A\to E_1\times E_2$ is dual to $p^*+q^*:E_1\times E_2\to A$ and that its
kernel $\Delta_A$ fits in the short exact sequence
$$0\to \Delta_E\to (E_1\times E_2)[2] \to \Delta_A\to 0$$ and, since
$(E_1\times E_2)[2]\simeq \Delta\times\Delta$, that
$\Delta_A\simeq\Delta$.

\begin{theorem}
\label{thrm:exact}
Let $E_1$ and $E_2$ be elliptic curves over a field $K$ of characteristic $0$
with isomorphic $m$-torsion $E_1[m]=\Delta=E_2[m]$. Let $A=(E_1\times
E_2)/\Delta_E$, where $\Delta_E$ is the anti-diagonal embedding of $\Delta$ in
$E_1[m]\times E_2[m]$
and $\mu_1: E_1\to H^1(K,\Delta)$ and $\mu_2: E_2\to H^1(K,\Delta)$. Then
\begin{itemize}
\item[(a)] $\mu_1\circ p_*=-\mu_2\circ q_*$
\item[(b)] The sequence
$$\xymatrix{
  E_1\times E_2(K)\ar[r]^(0.6){p^*+q^*}&
  A(K)\ar[r]^(0.4){\mu_1p_*}_(0.4){\mu_2q_*}&
  H^1(K,\Delta)}
$$
is exact.
\item[(c)] Let $\mu_1+\mu_2:E_1(K)\times E_2(K)\to H^1(K,\Delta)$ denote the map
$(e_1,e_2)\to \mu_1(e_1)+\mu_2(e_2)$. The sequence
$$\xymatrix{
  A(K)\ar[r]^(0.4){p_*\times q_*}&
  E_1\times E_2(K)\ar[r]^{\mu_1+\mu_2}&
  H^1(K,\Delta)}
$$
is exact.
\end{itemize}
\end{theorem}
\begin{proof}
(a): Let $a\in A(K)$ and put $\delta_1=\mu_1 p_*(a)$ and $\delta_2=\mu_2
q_*(a)$. By definition of $\mu_1$, there is a point $e_1\in E_1(\Kbar)$ such
that $me_1=p_*(a)$ and for any $\sigma\in \Gal(K)$, we have
$\act{\sigma}{e_1}=e_1+\delta_1(\sigma)$. Similarly, there is a point $e_2\in
E_2(\Kbar)$ with $me_2=q_*(a)$ and  $\act{\sigma}{e_2}=e_2+\delta_2(\sigma)$.

Since $p_*\circ(p^*\times q^*)=m|_{E_1}$ and $q_*\circ(p^*\times q^*)=m|_{E_2}$,
we see that $p^*(e_1)+q^*(e_2)-a=T\in\Delta_A(\Kbar)$. It follows that
$\act{\sigma}{T}-T=p^*\delta_1(\sigma)+q^*\delta_2(\sigma)$. Using the
identifications $E_1[m]=\Delta=E_2[m]$ via $p^*$ and $q^*$, we see that
$\delta_1+\delta_2$ is a 1-coboundary .

(b): It is immediate that $\mu_1 p_*\circ(p^*+q^*)=\mu_1\circ m|_{E_1}=0$.
On the
other hand, if $a\in A(K)$ with $\mu_1p_*(a)=0$ then, we can apply the same
construction as in (a), but now we can choose $e_1\in E_1(K)$ and $e_2\in
E_2(K)$. It follows that $p^*(e_1)+q^*(e_2)-a=T\in\Delta_A(K)$, so
$(p^*+q^*)(e_1-T,e_2)=a$.

(c): Using (a), we see that $(\mu_1+\mu_2)\circ(p_*\times q_*)=0$. To obtain
exactness, consider $e_1\in
E_1(K)$ and $e_2\in E_2(K)$ with $\mu_1(e_1)=-\mu_2(e_2)$. There are points
$e_3\in E_1(\Kbar)$ and $e_4\in E_2(\Kbar)$ such that $m e_3=e_1$, $m e_4=e_2$
and there is a point $T\in \Delta(\Kbar)$ such that
$p^*(\act{\sigma}{e_3}-e_3)+q^*(\act{\sigma}{e_4}-e_4)=\act{\sigma}{T}-T$
for all
$\sigma\in\Gal(K)$.
It follows that $a=p^*(e_3)+q^*(e_4)-T$
satisfies $\act{\sigma}{a}-a=0$, so $a\in A(K)$
and $p_*(a)=e_1$ and $q_*(a)=e_2$.
\end{proof}

\begin{corollary}
\label{cor:selmer}
Let $E_1$ and $E_2$ be elliptic curves over a number field $K$ with
$\Delta:=E_1[m]\simeq E_2[m]$. Let $p^*+q^*:E_1\times E_2\to A=(E_1\times
E_2)/\Delta$ and let $p_*\times q_*:A\to E_1\times E_2$ be the dual isogeny.
Then
\begin{itemize}
\item[(a)] $S^{(p_*\times q_*})(A/K)=S^{(m)}(E_1/K)+S^{(m)}(E_2/K)$,
\item[(b)] $S^{(p^*+q^*)}((E_1\times E_2)/K)=S^{(m)}(E_1/K)\cap S^{(m)}(E_2/K)$.
\end{itemize}
\end{corollary}

\begin{proof}
For all places $p$ of $K$, apply Theorem~\ref{thrm:exact} to the completion
$K_p$ of $K$ and use the definition of the Selmer group.
\end{proof}

From Corollary~\ref{cor:selmer}, and the fact that $m|_{E_1}$, $m|_{E_2}$,
$p^*+q^*$ and $p_*\times q_*$ all have kernels isomorphic to $\Delta$, it
follows that, for $m$ prime,
$$\rk S^{(p_*\times q_*)}(A/K)+\rk S^{(p^*+q^*)}((E_1\times E_2)/K)=
\rk S^{(m)}(E_1/K)+\rk S^{(m)}(E_2/K).$$
Furthermore, it is immediate that
$$S^{(m)}((E_1\times E_2)/K)=S^{(m)}(E_1/K)\times S^{(m)}(E_2/K),$$
so descents via $m|_{E_1\times E_2}$, the dual isogenies $m|_{E_1}\times
1|_{E_2}$ and $1|_{E_1}\times m|_{E_2}$, and the dual isogenies $p^*+q^*$ and
$p_*\times q_*$ all give rise to the same bound on $\rk (E_1\times E_2)(K)=
\rk E_1(K)+\rk E_2(K)$.

The Selmer-group $S^{(m)}(A/K)$ may yield different information. The map
$p_*\times q_*:A[m]\to \Delta$ induces a homomorphism $N$ that fits in the
following commutative diagram with exact rows.
$$\xymatrix{
A(K) \ar[r]^{m} \ar[d]^{p_*\times q_*} &
  A(K) \ar@{=}[d] \ar[r]^(0.4){\nu} &
  H^1(K,A[m]) \ar[d]^N\\
E_1\times E_2(K) \ar[r]^(0.6){p^*+q^*} &
  A(K) \ar[r]^(0.45){\mu_1 p_*} &
  H^1(K,\Delta)
}$$
This allows us to obtain a sharper bound on $\#A(K)/(p^*E_1(K)+q^*E_2(K))$.
Furthermore, using that $\mu_1 p_*(A(K))=N\nu(A(K))=\mu_2q_*(A(K))$, we find
that $N\nu(A(K))=\mu_1(E_1(K))=\mu_2(E_2(K))$. Consequently,
\begin{lemma}
\label{lemma:selint}
Let $A$, $E_1$ and $E_2$ be as defined above over a number field $K$, with $N:
H^1(K,A[m])\to H^1(K,\Delta)$. Then
$$\mu_1(E_1(K))\cap\mu_2(E_2(K))\subset S^{(m)}(E_1/K)\cap S^{(m)}(E_2/K) \cap
  NS^{(m)}(A/K)$$
\end{lemma}

Now we analyse when Lemma~\ref{lemma:selint} may provide a strictly sharper
bound than $S^{(p^*+q^*)}((E_1\times E_2)/K)$. This only applies if
$\Sha(E_1/K)[m]$ or $\Sha(E_2/K)[m]$ is non-trivial. 
Suppose that $\delta\in H^1(K,\Delta)$ represents a cocycle
with a non-trivial image under $H^1(K,\Delta)\to H^1(K,E_1)$, but which is
trivial in $H^1(K_p,E_1)$ for any place $p$ of $K$.
If we combine the Galois-cohomology of 
$E_1\stackrel{m}{\to}E_1$,
$E_2\stackrel{m}{\to}E_2$ and
$0\to E_1 \stackrel{p^*}{\to} A \stackrel{q_*}{\to} E_2 \to 0$, we obtain the
following commutative diagram with exact rows and columns.
$$\xymatrix{
&& E_2(K) \ar[r]^{q^*} \ar[d]_m &
  A(K) \ar[d]^{q_*} \\
&& E_2(K) \ar@{=}[r] \ar[d]_{\mu_2} &
  E_2(K) \ar[d]\\
E_1(K) \ar[r]^m \ar[d]_{p^*} &
  E_1(K)  \ar[r]^{\mu_1} \ar@{=}[d] &
  H^1(K,\Delta) \ar[r] \ar[d] &
  H^1(K,E_1) \ar[d]\\
A(K) \ar[r]^{p_*} &
  E_1(K) \ar[r] &
  H^1(K,E_2) \ar[r] &
  H^1(K,A)
}$$
If $\delta$ does not vanish in $H^1(K,A)$, it leads to a non-trivial element in
$\Sha(A/K)[p_*\times q_*]\subset\Sha(A/K)[m]$.

The space $H^1(K,E_1)$ is isomorphic to the set of principally homogeneous
spaces of $E_1$ over $K$ (see \cite[Theorem~X.3.6]{sil:AEC1}). As such
the map $E_2(K)\to H^1(K,E_1)$ corresponds to
$q_*^{-1}$, i.e., take the fibre of $q_*: A\to E_2$ over a rational point on
$E_2$. Following \cite{cremaz:vis}, we say that an element of $H^1(K,E_1)$
that occurs as such a fibre, is \emph{visualised} in $A$. We formulate this
observation as a corollary.

\begin{corollary} A cocycle $\xi\in H^1(K,E_1)$ vanishes in $H^1(K,A)$ precisely
if $\xi$ is in the image of $E_2(K)$.  
\end{corollary}

We see that a necessary condition for Lemma~\ref{lemma:selint} to yield an
improved rank bound, is that some non-trivial elements of
$\Sha(E_1/K)[m]\subset H^1(K,E_1)$ occur as fibres of $A$ over $E_2(K)$.

\section{Quadratic Weil-restrictions of elliptic curves}
\label{sec:quadext}

In this section, we interpret Theorem~\ref{thrm:quadext} in terms of the
construction explained in Section~\ref{sec:cohom} and we give a proof.
Informally, the idea is the following.
Recall from Section~\ref{sec:selell} that an element $\delta\in H^1(K,E[2])$ is
in the image of $E(K)/2E(K)$, and therefore maps to $0$ in $H^1(K,E)$, precisely
if a curve $T_\delta$ has a $K$-rational point. One obvious way to force a
rational point on $T_\delta$ is by extension of the base field. As is noted in
that same section, if $\delta\in S^{(2)}(E/K)$, then $T_\delta$ has a model
$u_3^2=R(t)$, where $R(t)$ is a quartic polynomial. For any value
$t_0\in K$ and $L:=K(\sqrt{R(t_0)})$,
the set $T_\delta(L)$ is non-empty. By taking the
\emph{Weil-restriction of scalars} $A=\Re_{L/K} E$, we find
ourselves in the situation of Section~\ref{sec:cohom} with $m=2$.
We make this explicit.

Let $K$ be a field of characteristic $0$. Consider the elliptic curves
$$\begin{array}{rrcl}
E_1:& y_1^2&=&x^3+a_2 x^2+a_4 x+ a_6=F(x)\\
E_2:& dy_2^2&=&x^3+a_2 x^2+a_4 x+ a_6\\
\end{array}$$
over $K$, where $d\in K^*$ is not a square. Let $L=K(\sqrt{d})$. As curves over
$L$, we have
$E_1\simeq E_2$ by $(x,y_1)=(x,\sqrt{d}y_2)$. We write $E(L)$ for $E_1(L)\simeq
E_2(L)$. We will use the model of $E_1$ for $E$, but the reader should note that
the construction is essentially symmetric in $E_1$ and $E_2$.

Obviously, $\Delta=E_1[2]\simeq E_2[2]$ as a $K$-Galois module. Applying
Theorem~\ref{thrm:H1E}, we find that the maps $E_i(K)/2E_i(K)\to H^1(K,\Delta)$
are induced by
$$\begin{array}{rrcl}
\mu_1: &E_1(K)&\to&M_K'\\
&(x,y_1)&\mapsto&(x-\theta)
\end{array}\quad
\begin{array}{rrcl}
\mu_2: &E_2(K)&\to&M_K'\\
&(x,y_2)&\mapsto&d(x-\theta)
\end{array}.$$

We consider the Weil-restriction $A=\Re_{L/K} E$ in the sense of
\cite[\S 7.6]{BLR:neron}. We suffice in describing $A(\Kbar)$ as a
$\Gal(\Kbar/K)$-module. We define $A(\Kbar)=E_1(\Kbar)\times E_1(\Kbar)$ as a
set, but with a twisted Galois-action.
For $\sigma\in\Gal(K)$ and $(e_1,e_2)\in A(\Kbar)$ we define
$$\act{\sigma}{(e_1,e_2)}=
\left\{
\begin{array}{ll}
(\act{\sigma}{e_1},\act{\sigma}{e_2})&
  \mbox{ if }\sigma\in\Gal(L)\subset\Gal(K)\\
(\act{\sigma}{e_2},\act{\sigma}{e_1})&
  \mbox{ otherwise.}
\end{array}
\right.$$

Let $\sigma\in \Gal(K)\setminus\Gal(L)$.
It is straightforward to check that $E(L)\simeq A(K)$ via
$e\mapsto(e,\act{\sigma}{e})$. We identify $E_2(\Kbar)$ with $E_1(\Kbar)$ via
$(x,y_1)=(x,\sqrt{d}y_2)$. Under this identification, we have
$$E_2(K)=\{e\in E_1(\Kbar):e=-\act{\sigma}{e}\mbox{ for
}\sigma\in\Gal(K)\setminus\Gal(L)\}.$$

As maps between Galois-modules, we obtain.
$$\begin{array}{rccc}
p^*:&E_1(\Kbar)&\to&A(\Kbar)\\
&e&\mapsto&(e,e)
\end{array}
\begin{array}{rccc}
q^*:&E_1(\Kbar)&\to&A(\Kbar)\\
&e&\mapsto&(e,-e)
\end{array}$$
$$\begin{array}{rccc}
p_*:&A(\Kbar)&\to&E_1(\Kbar)\\
&(e_1,e_2)&\mapsto&e_1+e_2
\end{array}
\begin{array}{rccc}
q_*:&A(\Kbar)&\to&E_1(\Kbar)\\
&(e_1,e_2)&\mapsto&e_1-e_2
\end{array}$$
It is straightforward to check that this places us in the situation of
Section~\ref{sec:cohom} with $m=2$.

By Shapiro's Lemma and Theorem~\ref{thrm:H1E},
we have $H^1(K,A[2])=H^1(L,E_1[2])=M_L'/M_L^*$. We write $\nu:A(K)/2A(K)\to
M_L'/M_L^*$ for the corresponding map, where we use $A(K)\simeq E(L)$
and the map $\res^K_L\mu_1: E_1(L)\to M'_L$.

The obvious map $N=N_{M_L/M_K}: M_L^*\to M_K^*$ completes the explicit
description of the ingredients of Corollary~\ref{cor:selmer} and
Lemma~\ref{lemma:selint}, as the following lemma shows.

\begin{lemma} With the definitions above,
$$\mu_1 p_*= N\nu = \mu_2 q_*$$
\end{lemma}
\begin{proof}
We prove the first equality. The second follows by symmetry.
Suppose $a\in A(K)$ and let $e\in E(L)$ be the corresponding point. Then $\mu_1
p_*(a)=\mu_1(e+\act{\sigma}{e})$, where $\sigma\in\Gal(K)$ acts non-trivially on
$L$. By Corollary~\ref{cor:munorm}, this equals $N\nu(a)$.
\end{proof}

\noindent\emph{Proof of Theorem~\ref{thrm:quadext}}: Let $\delta\in S^{(2)}$ be
a cocycle representing $\xi$. Recall from Section~\ref{sec:selell} that
$L_\delta$ is a curve of genus $0$ with points everywhere locally. Consequently,
$L_\delta$ is isomorphic to $\PP^1$ and has infinitely many rational points.
Each point lifts to a quadratic point on the genus $1$ curve $T_\delta$. Thus
$T_\delta(L)$ is non-empty for some quadratic extension of $K$. Since the number
of $L$-rational points of $T_\delta$ of height bounded by $B$ is at most
logarithmic in $B$ and the number of $K$-rational points on $L$ of height
bounded by $B$ is polynomial in $B$, we see that this construction leads to
infinitely many distinct quadratic extensions $L$ of $K$.

Note that, given $t\in L_\delta(K)$, we can take $d=F(x(t))$. Then, $E^{(d)}:
dy_2^2=F(x)$ has a rational point over $x(t)$. Therefore, the fibre product
$T_\delta^{(d)}=E^{(d)}\times_{\PP^1} L_\delta$ has a rational point over $x(t)$ as
well. Since, over $L=K(\sqrt{d})$, the curves $E$ and $E^{(d)}$ are isomorphic,
we see that $T_\delta(L)\simeq T_\delta^{(d)}(L)$ is non-empty.

\section{Bielliptic curves of genus $2$}
\label{sec:biell}

The previous section shows that for an elliptic curve $E$ over a number field
$K$, any element of $\Sha(E/K)[2]$ can be visualised in infinitely many distinct
Abelian surfaces. The Abelian surfaces constructed there are all isomorphic to
$E\times E$ over $\Kbar$, however. In this section we prove
Theorem~\ref{thrm:biell}, i.e., that we can insist on non-isomorphic surfaces
over $\Kbar$. 

Let $K$ be a field of characteristic $0$ and consider
$$\begin{array}{rrcl}
E_1:& y_1^2&=&x^3+a_2x^2+a_4x+a_6=F(x)\\
L_0:&dy_0^2&=&x-a\\
E_2:&dy_2^2&=&(x-a)F(x).
\end{array}$$
We consider these curves to be covers of the same $\PP^1$ corresponding to $x$.
The fibre-product $C=E_1\times_{\PP^1} L_0=E_2\times_{\PP^1} L_0$ is a curve of
genus $2$ with model
$$C: y_1^2=F(dy_0^2+a).$$
The curve $C$ is a double cover of both $E_1$ and $E_2$ by
$$\begin{array}{rcccl}
p:& C&\to&E_1\\
&(y_0,y_1)&\mapsto&(x,y_1)&=(dy_0^2+a,y_1)\\[0.5em]
q:& C&\to&E_2\\
&(y_0,y_1)&\mapsto(x,y_2)&=(dy_0^2+a,y_0y_1)
\end{array}$$
We write $\infty\in L_0$ and $\infty\in E_1$ for the unique points on $L_0$ and
$E_1$ above $\infty\in\PP^1$. Choosing an arbitrary but fixed labelling, we
write $\infty^+,\infty^-\in C$ for the two points on $C$ above $\infty\in L_0$.
The points $\infty^+$ and $\infty^-$ are
rational or quadratic conjugate over $K$, depending on $d$ being a square. The
unique point on $E_1$ above $\infty\in \PP^1$ is denoted by $\infty$ as well.

For a curve $C$ over $K$, we represent an element of $A:=\Jac_C$ by the
corresponding divisor class $[\sum_{i=1}^r n_i P_i]$, where $P_i\in C(\Kbar)$
and $n_i\in \ZZ$, with $n_1+\cdots+n_r=0$. Since we only consider curves of
genus $1$ and $2$, we have that any point in $\Jac_C(K)$ can be represented by
the class of a $K$-rational divisor, i.e. a $\Gal(K)$-stable formal linear
combination of points in $C(\Kbar)$.

We identify $E_1$ with $\Jac_{E_1}$ via $e_1\mapsto[e_1-\infty]$ and $E_2$ with
$\Jac_{E_2}$ via $e_2\mapsto[e_2-(a,0)]$. The map $p$ gives rise to maps
between Abelian varieties
$$\begin{array}{rccc}
p_*:&\Jac_C&\to&E_1\\
&\left[\sum_{i=1}^r n_i P_i\right]&\mapsto&\sum_{i=1}^r n_i p(P_i)\\[0.5em]
p^*:&E_1&\to&\Jac_C\\
&e&\mapsto&\left[\sum_{P\in p^{-1}\{e\}} P -\infty^+ - \infty^-\right]
\end{array}$$
and similarly for $q$. It is straightforward to check that these maps satisfy
the assumptions in Section~\ref{sec:cohom} and we adopt the notation from that
section. 

Let $\theta\in M_K=K[x]/F(x)$ be such that $F(\theta)=0$. Using
Theorem~\ref{thrm:H1E} we find
$$\begin{array}{rcccl}
\mu_1:&E_1(K)/2E_1(K)&\to&M_K'/M_K^{*2}\\
 &(x,y_1)&\mapsto&(x-\theta)&\mbox{ if }y_1\neq 0
\end{array}$$
and, by transforming $E_2$ into Weierstrass form,
$$\begin{array}{rcccl}
\mu_2:&E_2(K)/2E_2(K)&\to&M_K'/M_K^{*2}\\
 &(x,y_1)&\mapsto&d(x-a)(a-\theta)F(a)(x-\theta)&\mbox{ if }y_2\neq 0\\
 & & &dF(a)(a-\theta)&\mbox{ if }x=\infty
\end{array}$$

Any point $a\in A(K)$ has a representative
$a=[(y_{0,1},y_{1,1})+(y_{0,2},y_{1,2})-\infty^+-\infty^-]$, where
$(y_{0,1},y_{1,1})$ and $(y_{0,2},y_{1,2})$ are either rational points of $C$ or
quadratic conjugates. In fact, this representative is unique if $a\neq 0$.
Let $\sM_K=K[y_0]/F(dy_0^2+a)$ and let $\Theta$ be the class of $X$ in $\sM_K$.
Similar to the map $\mu$ we defined in Theorem~\ref{thrm:H1E}, we define
$$\begin{array}{cccc}
\tilde{\nu}:&A(K)&\to&\sM_K'/K^*\sM_K^{*2}\\
&a&\mapsto&(y_{0,1}-\Theta)(y_{0,2}-\Theta)
\end{array}$$
Unfortunately, the map $\tilde{\nu}$ does not completely correspond to
$\nu: A(K)\to H^1(K,A[2])$. The kernel may be a little bigger.
\begin{lemma}[\cite{prolegom}, \cite{stoll:descent}]
With the definitions above, there is a map $H^1(K,A[2])\to \sM_K'/K^*\sM_K^{*2}$
that makes the following diagram with exact row commutative.
$$\xymatrix{
A(K)\ar[r]^2&A(K)\ar[r]^\nu\ar[dr]^{\tilde{\nu}}&H^1(K,A[2])\ar[d]\\
&&\sM_K'/K^*\sM_K^{*2}
}$$
\end{lemma}
For $K$ a number field, we define the \emph{fake Selmer group}
$\tilde{S}^{(2)}(A/K)$ to be the image of
$S^{(2)}(A/K)$ under this map $H^1(K,A[2])\to \sM_K'/K^*\sM_K^{*2}$. This group
is effectively computable (see \cite{stoll:descent}).
Fortunately, the map $N:H^1(K,A[2])\to H^1(K,\Delta)$ factors through
$\sM_K'/K^*\sM_K^{*2}$, so for applying Lemma~\ref{lemma:selint}, knowing the
fake Selmer group suffices.
\begin{lemma}
With the definitions above, the following identity holds.
$$\mu_1 p_*=N_{\sM_K/M_K}\tilde{\nu}= \mu_2 q_*$$
\end{lemma}
\begin{proof} We prove the first equality. The second follows by symmetry.
Suppose $a=[(y_{0,1},y_{1,1})+(y_{0,2},y_{1,2})-\infty^+-\infty^-]\in A(K)$.
We have $p_*(a)=(x_1,y_{1,1})+(x_2,y_{1,2})\in E_1(K)$, where $x_1=dy_{0,1}^2+a$
and $x_2=dy_{0,2}^2+a$. Generically, by Corollary~\ref{cor:munorm}, we have
$$\mu_1 p_*(a)=(x_1-\theta)(x_2-\theta)=N_{\sM_K/M_K}\tilde{\nu}(a).$$
The special points where the middle expression is singular, can be tested
separately.
\end{proof}

\noindent\emph{Proof of Theorem\ref{thrm:biell}}:
Let $\delta\in H^1(K,E_1[2])$ be a cocycle mapping to $\xi$. Since $T_\delta$ as
defined in Section~\ref{sec:selell} has points everywhere locally, the subcover
$L_\delta\simeq \PP^1$. Take a rational point $t_0\in L_\delta$.

For some $a\in K$ with $F(a)\neq 0$ and $d\in K^*$, we consider
$$E_2:dy_2^2=(x-a)F(x).$$
The cover of $E_2$ corresponding to $\delta\in H^1(K,E_2[2])\simeq
H^1(K,E_1[2])$ is $L_\delta\times_{\PP^1}E_2$. If we choose
$d=(x(t_0)-a)F(x(t_0)),$
then this cover has a rational point above $x(t_0)$. We see that we are
completely free in choosing $a$, which gives us infinitely many distinct curves
$E_2$ over $\Kbar$ with the desired property.

\noindent\textbf{Remark:} By letting $a\to\infty$, we get the construction from
Section~\ref{sec:quadext} as a limit of the construction in this section.

\section{Algorithmic considerations}

The fact that the constructions presented in Sections~\ref{sec:quadext} and
\ref{sec:biell} are completely explicit and that the Selmer groups in
Lemma~\ref{lemma:selint} are effectively and often practically computable,
suggests a probabilistic algorithm to prove non-triviality of $\Sha(E/K)[2]$ for
an elliptic curve $E$ over a number field $K$.

Let $\delta\in M_K^*$ represent an element of $S^{(2)}(E/K)$ for which we
suspect that $T_\delta(K)$ is empty, i.e., $\delta$ represents a non-trivial
element in $\Sha(E/K)[2]$.

\begin{enumerate}
\item Choose $t_0\in L_\delta(K)$ and compute $x_0=x(t_0)$.
\item Take $E_2: dy_2=F(x)$ or $E_2:dy_2^2=(x-a)F(x)$ such that $E_2$ has a
rational point above $x_0$.
\item Test if $\delta\in S^{(2)}(E/K)\cap S^{(2)}(E_2/K)\cap N S^{(2)}(A/K)$.
\item If this is the case, go to (1). Otherwise, $\delta$ indeed represents a
non-trivial element in $\Sha(E/K)[2]$.
\end{enumerate}

Of course, if $T_\delta$ has a rational point then this algorithm will not
terminate. Unfortunately, otherwise we cannot prove it will either. We cannot
exclude that $\delta \notin NS^{(2)}(A/K)$.
In practice, the procedure appears to have a rather high success rate, though.

Of course, if $\delta$ represents a non-trivial element of $2\Sha(E/K)[4]$, then
the algorithm will not terminate either. The image of $\delta$ in $\Sha(A/K)[2]$
will indeed be trivial, but a class which corresponds to $\frac{1}{2}\delta$
will just map
to a non-trivial element in $\Sha(A/K)[2]$. However, if we have $A=\Re_{L/K} E$
as in Section~\ref{sec:quadext} then $\Sha(A/K)[2]\simeq \Sha(E/L)[2]$. We can
use the same algorithm to prove that the image of that class in $\Sha(A/K)$
indeed is non-trivial. See Example~4 in Section~\ref{sec:examples}

Theoretically, this procedure is weaker than attempting a sequence of $2$-power
descents, as suggested in for instance \cite[Proposition~X.4.12]{sil:AEC1}.
The latter is guaranteed to work if $\Sha(E/K)[2^r]=0$ for some $r$. The
procedure above needs that and then still is only probabilistic. See for
instance Example~2 in Section~\ref{sec:examples}

For practical applications, this approach has the benefit that one computes
$2$-Selmer groups of an elliptic curve over a tower of quadratic extensions,
rather than higher Selmer groups over a constant base field.
See for instance \cite{mss:4desc} for $4$-descents and
\cite{cassels:seconddesc} second descents, that yield the same information as
$4$-descents. See also \cite{siksma:compsel} for a complexity analysis of
$2$-descents and \cite{dses:pselm} and \cite{schasto:pdescent} for odd
$p$-descents.

The computation of a $2$-Selmer group of an elliptic curve over an arbitrary
number field, although for higher degree fields prohibitively laborious, is
effectively implemented in MAGMA \cite{magma} and KASH \cite{kant} using the
routines for computing $S$-unit groups originally by Hess \cite{hess:diplom}.
See \cite{bruin:visshamag} and \cite{bruin:algae}. Such an implementation is
presently missing for higher $2$-power Selmer groups.

\section{Examples}
\label{sec:examples}
As an illustration of the practicality of the methods described in this article,
we give some examples. We only give a brief outline of the computations, to show
the reader what phenomena occur. For more detailed results, we refer the reader
to \cite{bruin:visshamag}, which contains a complete transcript of the MAGMA
\cite{magma} sessions leading to the results presented here. If the reader has
access to MAGMA, the software to repeat the computations is also included.

\noindent\textbf{Example 1. Mutual visualisation.} Consider the curve from
Proposition~\ref{prop:schasto}. A computation shows
$$S^{(2)}(E_1/\QQ)=\langle -\theta, 1-\theta, 1-2\theta, 1-8\theta\rangle
\mbox{ and }
\mu_1\langle (0,1),(1,1)\rangle=\langle -\theta,1-\theta\rangle.$$
We suspect that $1-2\theta$ and $1-8\theta$ represent non-trivial elements in
$\Sha(E_1/\QQ)[2]$. If we follow the procedure above, we write down the curve
$L_{1-2\theta}$,find a rational point $t_0$ on it and compute $d=F(x(t_0))$.
From computational considerations, it is desirable that $d$ has a small
square-free part. Instead of trying different $t_0$ to obtain such a $d$, we
enumerate small $x$-coordinates and check if the point lifts to $L_{1-2\theta}$.

We find that $L_{1-2\theta}$ has a rational point above $x=\frac{1}{2}$ and that
$L_{1-8\theta}$ has a rational point above $x=\frac{1}{8}$. Furthermore,
$F(\frac{1}{2})=2(\frac{7}{4})^2$ and $F(\frac{1}{8})=2(\frac{41}{32})^2$, so if
we take $d=2$ in the construction of Section~\ref{sec:quadext}, both elements
will be visualised in the Abelian surface $A$. In fact, we find
$$S^{(2)}(E_1/\QQ(\sqrt{2}))=\langle -\theta, 1-\theta, 1-2\theta,
1-8\theta\rangle,$$
so
$$NS^{(2)}(E_1/\QQ(\sqrt{2}))=\langle 1\rangle.$$

By Lemma~\ref{lemma:selint}, we see that
$\mu_1(E_1(\QQ))\cap\mu_2(E_2(\QQ))=\langle 1\rangle$. Since $\langle -\theta,
1-\theta\rangle\subset \mu_1(E_1(\QQ))$ and $\langle 1-2\theta,
1-8\theta\rangle\subset \mu_2(E_2(\QQ))$, we see that equality must hold.

\noindent\textbf{Example 2. Failed visualisation.} In the previous example we
were lucky that both elements of $\Sha(E_1/\QQ)[2]$ are visualised
simultaneously. This need not happen. For the same curve as above, consider
$L_{(1-\theta)(1-2\theta)}$. It has a rational point above $\frac{-4}{11}$.
Since $F(-\frac{4}{11})=-11(\frac{113}{121})^2$, it follows that
$(1-\theta)(1-2\theta)$ is visualised in the surface constructed according to
Section~\ref{sec:quadext} with $d=-11$. As it turns out,
$$S^{(2)}(E_2/\QQ)=\mu_2E_2(\QQ)=\langle (1-\theta)(1-2\theta)\rangle.$$
and
$$NS^{(2)}(E_1/\QQ(\sqrt{-11}))=\langle (1-\theta)(1-2\theta)\rangle.$$
While we know we have visualised an non-trivial element from $\Sha(E_1/\QQ)[2]$
from the argument given above, we cannot conclude this from the data computed
with $d=-11$.

\noindent\textbf{Example 3. Visualisation using bielliptic curves of genus $2$.}
The construction from Section~\ref{sec:biell} offers two degrees of freedom,
giving enhanced flexibility. We again consider the curve $E_1:
y^2=F(x)=x^3-22 x^2+21 x+1$. 

We pick two elements from $S^{(2)}(E/\QQ)$, say $\delta_1=1-2\theta$ and
$\delta_2=-\theta(1-2\theta)(1-8\theta)$. We pick $x_1=1/2$ below
$L_{\delta_1}(\QQ)$ and $x_2=9/10$ below $L_{\delta_2}(\QQ)$. 
We use the construction from Section~\ref{sec:biell} to obtain an Abelian
surface $A=\Jac_C$ in which both $\delta_1$ and $\delta_2$ are visualised. To
this end, we determine $a,d\in\QQ$ with $d\neq 0$ and $F(a)\neq 0$
such that both $d(x_1-a)F(x_1)$ and $d(x_2-a)F(x_2)$ are squares. A priori such
$a,d$ might not exist. In that case, one could choose other $x_1$ and $x_2$ and
try again. While there is no guarantee that this procedure yields success after
finitely many steps, in practice it seems rather easy to find $x_1,x_2,a,d$
with the desired properties. In our case $a=1$, $d=-1$ works. Thus $\delta_1$
and $\delta_2$ are visualised in the the Abelian
surface $A=\Jac_C$, where $C:y_1^2=-y_0^6-19y_0^4+20y_0^2+1$. We have
$$S^{(2)}(E_1/\QQ)\simeq(\ZZ/2\ZZ)^4\mbox{ and }
  S^{(2)}(E_2/\QQ)\simeq(\ZZ/2\ZZ)^3,$$
with
$$S^{(2)}(E_1/\QQ)\cap S^{(2)}(E_2/\QQ)\simeq(\ZZ/2\ZZ)^2.$$
Using Stoll's routines
\cite{stoll:descent} in MAGMA \cite{magma}, 
we can compute the fake Selmer group
$\tilde{S}^{(2)}(\Jac_C/\QQ)\simeq(\ZZ/2\ZZ)^4$
and
$$NS^{(2)}(\Jac_C/\QQ)=\langle 21-\theta\rangle.$$
Writing
$$E_2:v^2=u^3+20u^2-19u-1,\mbox{ where }(u,v)=(1/y_0^2,y_1/y_0^3),$$
we find
$$\mu_2\langle (1,1),(2,7),(5,23)\rangle=
  \langle -\theta(1-2\theta),-\theta(1-8\theta),21-\theta\rangle.$$
By Lemma~\ref{lemma:selint}, we see that $1-2\theta$ and $1-8\theta$ in
$S^{(2)}(E_1/\QQ)$ are not in $\mu_1(E_1(\QQ))$.

\noindent\textbf{Example 4. Visualisation of $4$-torsion.}
We consider the curve
$$E_1: y_1^2=x^3-162x^2+x.$$
This is a model of the curve 1640G3 in \cite{cremona}. Using a $4$-isogeny,
one can show that $E_1(\QQ)=\ZZ/2\ZZ$ and that $\Sha(E_1/\QQ)[4]=(\ZZ/4\ZZ)^2$.
We use this curve to illustrate Corollary~\ref{cor:twopowvis}, which also
applies to curves without any $2$-power isogenies.

Let $E_2:-y_2^2=x^3-162x^2+x$ and $K=\QQ(i)=\QQ(\sqrt{-1})$. We consider a
subgroup
$\langle g_1,\ldots,g_8\rangle\subset H^1(\QQ,E[2])$, which contains both
$S^{(2)}(E_1/K)$ and $S^{(2)}(E_2/K)$, corresponding to generators chosen by 
\cite{bruin:visshamag} based on \cite{magma}. Since we merely want to
illustrate Corollary~\ref{cor:twopowvis}, we refer the reader to
\cite{bruin:visshamag} for details and the definitions of the $g_i$.
We find
$$\begin{array}{lcl}
S^{(2)}(E_1/\QQ)=\langle g_5,g_6,g_8\rangle&;&\mu_1(0,0)=g_5+g_6\\
S^{(2)}(E_2/\QQ)=\langle g_6,g_6+g_8\rangle&=&
                \mu_2\langle(0,0),(-9/4,231/8)\rangle.
\end{array}$$
In particular, the group $\langle g_6,g_8\rangle\subset S^{(2)}(E_1/\QQ)$ is
visualised in $\Re_{K/\QQ}E_1$. However, we already know that this group
represents $2\Sha(E_1/\QQ)[4]$. Indeed, writing $\langle h_1,\ldots, h_7\rangle\subset
H^1(K,E[2])$, we find
$$S^{(2)}(E_1/K)=\langle h_3,h_4,h_5+h_6,h_7\rangle;\;
\res^\QQ_K\mu_1\langle(0,0),(-9/4,231/8i)\rangle=\langle h_4,h_3\rangle$$
and
$$N_{K/\QQ} S^{(2)}E_1/K)=\langle g_6,g_8\rangle.$$
Of course, $N_{K/\QQ} h_3=N_{K/\QQ} h_4=0$.

As suggested in the proof of Corollary~\ref{cor:twopowvis}, We pick an element
$h=h_4+h_5+h_6\in S^{(2)}(E_1/K)$ such that $N_{K/\QQ} h = g_8$ is suspected
to represent a non-trivial element in $\Sha(E/\QQ)[2]$. We choose
$$E_3: (i-2)y_3^2=x^3-162x^2+x,$$
such that under the map $\mu_3:E_3(K)\to H^1(K,E[2])$ we have
$$\mu_3(0,0)=h\in S^{(2)}(E_3/K).$$
Therefore, putting $L=K(\sqrt{i-2})$, the element $h\in S^{(2)}(E_1/K)$ is
visualised in $\Re_{L/K} E_1$. Indeed, we find that
$$N_{L/K} S^{(2)}(E_1/L)=\langle h_4,h_3+h_5+h_6\rangle.$$
Using Lemma~\ref{lemma:selint}, we find that
$\res^\QQ_K\mu_1(E_1(K))\cap\mu_3(E_3)\subset N_{L/K}S^{(2)}(E_1/L)$.
In particular, it follows that
$$h_5+h_6+\langle h_3,h_4\rangle\cap\res^\QQ_K\mu_1(E_1(K))=\emptyset.$$
Since
$$h_5+h_6+\langle h_3,h_4\rangle=\{h\in S^{(2)}(E_1/K): N_{K/\QQ}h=g_8\},$$
we see that $g_8\notin\mu_1(E_1(\QQ))$. It follows that $g_8$ indeed represents
a non-trivial element in $\Sha(E_1/\QQ)[2]$. We could give a similar argument
for $g_6$.

\providecommand{\bysame}{\leavevmode\hbox to3em{\hrulefill}\thinspace}

\end{document}